\documentclass[12pt,reqno]{amsart}
\usepackage{amsmath,amsthm,amsopn}
\usepackage{enumerate}
\usepackage[a4paper, margin={2.5cm}]{geometry}
\usepackage{mathrsfs}

\newtheorem{theorem}{Theorem}

\theoremstyle{definition}

\theoremstyle{remark}
\newtheorem{remark}{Remark}
\theoremstyle{problem}
\newtheorem{problem}{Problem}

\newcommand{\R}{\mathbb R}
\newcommand{\N}{\mathbb N}
\newcommand{\e}{\mathrm e}
\newcommand{\set}[1]{\left\{#1\right\}}

\newcommand{\yk}[1][k]{Y^{#1}}

\renewcommand{\P}[1][]{\mathbb P^{#1}}
\newcommand{\E}[1][]{\mathbb E^{#1}}

\newcommand{\bone}{\mathbf 1}
\newcommand{\bfd}{\rho}
\newcommand{\eF}{F }
\newcommand{\prt}{\partial}
\newcommand{\eps}{\varepsilon}
\newcommand{\wh}{\widehat}

\DeclareMathOperator{\dist}{dist}

\DeclareMathOperator{\bessel}{Bes}
\DeclareMathOperator{\besselsq}{BESQ}

\newcommand{\bes}[2][\nu]{\bessel^{#1}(#2)}
\newcommand{\besq}[2][\nu]{\besselsq^{#1}(#2)}
\newcommand{\X}{\mathbf X}
\newcommand{\x}{\mathbf x}
\newcommand{\Y}{\mathbf Y}

\newcommand{\setof}[2]{\bc{{#1}\,:\,{#2}}}
\newcommand{\br}[1]{\left( #1\right)}

\newcommand{\bc}[1]{\left\{ #1\right\}}
\newcommand{\halfint}[1]{\left[ #1 \right)}

\newcommand{\abs}[1]{\left\lvert #1 \right\rvert}

\numberwithin{equation}{section}

\begin{document}

\title{Extinction of Fleming-Viot-type particle systems with strong drift}

\author{
{\bf Mariusz Bieniek},  {\bf Krzysztof Burdzy}  \ and \ {\bf Soumik Pal} }

\address{MB: Instytut Matematyki, Uniwersytet Marii Sk\l odowskiej-Curie, 20-031
Lublin, Poland}
\address{KB, SP: Department of Mathematics, Box 354350,
University of Washington, Seattle, WA 98195, USA}

\email{mariusz.bieniek@poczta.umcs.lublin.pl}
\email{burdzy@math.washington.edu}
\email{soumik@math.washington.edu}

\thanks{Research supported in part by NSF Grants DMS-0906743, DMS-1007563, and by grant N N201
397137, MNiSW, Poland.}

\begin{abstract}
We consider a Fleming-Viot-type particle system consisting of independently moving particles that are killed on the boundary of a domain. At the time of death of a particle, another particle branches. If there are only two particles and the underlying motion is a Bessel process on $(0,\infty)$, both particles converge to 0 at a finite time if and only if the dimension of the Bessel process is less than 0. If the underlying diffusion is Brownian motion with a drift stronger than (but arbitrarily close to, in a suitable sense) the drift of a Bessel process, all particles converge to 0 at a finite time, for any number of particles. 
\end{abstract}


\maketitle

\section{Introduction}\label{sec:intro}

Our paper is motivated by an open problem concerning extinction in a finite time of a branching particle system. We prove two results that are related to the original problem and might shed some light on the still unanswered question. 

The following Fleming-Viot-type particle system was studied in
\cite{BurdzyMarch00}. Consider an open bounded set $D\subset\R^d$ and an integer $N\geq
2$. Let $\X_t=(X_t^1,\dotsc,X_t^N)$ be a process with values in $D^N$ defined as follows.
Let $\X_0=(x^1,\dotsc,x^N)\in D^N$. Then the processes $X_t^1,\dotsc,X_t^N$ evolve as
independent Brownian motions until the time $\tau_1$ when one of them, say, $X^j$ hits the
boundary of $D$. At this time one of the remaining particles is chosen uniformly, say,
$X^k$, and the process $X^j$ jumps at time $\tau_1$ to $X^k_{\tau_1}$. The processes
$X_t^1,\dotsc,X_t^N$ continue evolving as independent Brownian motions after time $\tau_1$
until the first time $\tau_2>\tau_1$ when one of them hits the boundary of $D$.  Again at
the time $\tau_2$ the particle which approaches the boundary jumps to the current location
of a particle chosen uniformly at random from amongst the ones strictly inside $D$.  The
subsequent evolution of $\X$ proceeds in the same way. 
We will say that $\X$ constructed above is \emph{driven} by Brownian motion. The main results in this paper are concerned with Fleming-Viot particle systems driven by other processes.

The above recipe defines the
process $\X_t$ only for $t < \tau_\infty$, where 
\begin{equation*}
  \tau_\infty = \lim_{k\to \infty} \tau_k. 
\end{equation*}
There is no natural way to define the process $\X_t$ for $t\geq \tau_\infty$, and,
therefore, it is of interest to investigate what conditions ensure that
$\tau_\infty=\infty$.  In Theorem 1.1 of \cite{BurdzyMarch00} the authors claim that in every domain
$D\subset\R^d$ and every $N\geq 2$, we have $\tau_\infty=\infty$, so the Fleming-Viot
process is always well-defined. However, the proof of Theorem 1.1 in \cite{BurdzyMarch00}
contains an error which is irreparable in the following sense. That proof is based on only two properties of Brownian
motion---the strong Markov property and the fact the the hitting time distribution of a
compact set has no atoms (assuming that the starting point lies outside the set). Hence,
if some version of that argument were true, it would apply to almost all non-trivial
examples of Markov processes with continuous time, and in particular to all diffusions.
However, in \cite{BieniekBurdzyFinch10}, the authors provided an example of a diffusion
$X$ on $D=(0,\infty)$ (a Bessel process with dimension
$\nu=-4$), such that $\tau_\infty<\infty$ for the Fleming-Viot process driven by this diffusion with $N=2$. 

It is not known whether Theorem 1.1 of \cite{BurdzyMarch00} is correct in full generality. It was proved in \cite{BieniekBurdzyFinch10,GK} that the theorem holds in domains which do not have thin channels.

\subsection{Main results}

We will prove two theorems. The first theorem is concerned with Bessel processes but it is motivated by the original model based on Brownian motion in an open bounded subset of $\R^d$. 
Recall that for any real $\nu$, a $\nu$-dimensional Bessel process on $(0,\infty)$ killed at 0 may be defined as a solution to the stochastic differential equation
\begin{align}\label{s30.3}
dX_t = dW_t + \frac{\nu-1}{2 X_t} dt,
\end{align}
where $W$ is the standard Brownian motion.
To make a link between Brownian motion in a domain and Bessel processes, we recall that there exists a regularized version $\bfd$ of the distance function (\cite[Theorem 2, p. 171]{Stein}). More precisely, there exist $0<c_1,c_2,c_3, c_4<\infty$ and a $C^\infty$ function $\bfd: D \to (0,\infty)$ with the following properties,
\begin{align*}
&c_1 \dist(x, \prt D)  \leq \bfd(x) \leq c_2 \dist(x, \prt D) ,\\
&\sup_{x\in D} |\nabla \bfd(x)| \leq c_3,  
\\
&\sup_{x\in D} \left| \bfd (x) \frac\prt{\prt x_i} \frac\prt{\prt x_m} \bfd(x)\right| \leq c_4   \qquad \hbox{for } 1\leq i,m\leq d.
\end{align*}
The above estimates and the It\^o formula show that if $B=(B^1, \dots, B^d)$ is a $d$-dimensional Brownian motion and $Z_t = \bfd(B_t)$ then 
\begin{align*}
d Z_t = \sum_{k=1}^ d a_k(Z_t) dB^k_t + \frac{b(Z_t)}{Z_t} dt,
\end{align*}
where the functions $a_k(\,\cdot\,)$ and $b(\,\cdot\,)$ are bounded. This shows that the dynamics of $Z$ resembles that of a Bessel process. 
Note that if $\tau_\infty < \infty$ for the Fleming-Viot process driven by Brownian motion in a domain $D$ then the distances of all particles to $\prt D$ go to 0 as $t\uparrow \tau_\infty$, by Lemma 5.2 of 
\cite{BieniekBurdzyFinch10}.
Hence, it is of some interest to see whether a Fleming-Viot process based on a Bessel process can become extinct in a finite time. We have a complete answer only for $N=2$, i.e., a two-particle process.

\begin{theorem}\label{thm:main1}
Let $\X$ be a Fleming-Viot process with $N$ particles on $(0,\infty)$ driven by Bessel process of dimension $\nu\in\R$. 

(i) If $N=2$
then $\tau_\infty<\infty$, a.s., if and only if $\nu< 0$.

(ii) 
If $N\nu \geq 2$ then $\tau_\infty=\infty$, a.s.
\end{theorem}

\bigskip
Our second main result is also motivated by some results presented in \cite{BurdzyMarch00}.
Several theorems in 
\cite{BurdzyMarch00,BieniekBurdzyFinch10} are concerned with limits when $N\to \infty$. To formulate rigorously any of these theorems it would suffice that $\tau_\infty = \infty$, a.s., for all sufficiently large $N$. In other words, it is not necessary to know whether $\tau_\infty = \infty$ for small values of $N$. One may wonder whether 
it is necessarily the case that $\tau_\infty = \infty$ for any Fleming-Viot-type process and sufficiently large $N$. Our next result shows that once the drift of the diffusion is slightly stronger than the drift of any Bessel process then $\tau_\infty < \infty$ for the Fleming-Viot process driven by this diffusion and {\it every} $N$.

Consider the following SDE for a diffusion on $(0,2]$,
\begin{equation}\label{eq:SDEX1}
 X_t=x_0+W_t-\int_0^t  \frac1{\beta X_t^{\beta -1}}\,ds-L_t,\quad t\leq T_0,
\end{equation}
where $x_0 \in (0,2]$, $\beta >2$, $W$ is Brownian motion, $T_0$ is the first hitting time of 0 by $X$, and $L_t$ is the local time of
$X$ at $2$, i.e., $L_t$ is a CAF of $X$ such that
\begin{equation*}
  \int_0^\infty \bone_{\set{X_s\ne 2}}dL_s=0,\quad\text{a.s.}
\end{equation*}
It is well known that \eqref{eq:SDEX1} has a unique pathwise solution $(X,L)$ (see, e.g.,
\cite{BassDEO}, Theorem I.12.1).
We will analyze a Fleming-Viot process on $(0,2]$ driven by the diffusion defined in \eqref{eq:SDEX1}. The role of the boundary is played by the point 0, and only this point.
In other words, the particles jump only when they approach 0.
Let $\P[\x]$ denote the distribution of 
the Fleming-Viot particle system starting from $\X_0 =\x$.

\begin{theorem}\label{thm:main2}
 Fix any $\beta >2$. For every $N\geq 2$, the $N$-particle Fleming-Viot process on $(0,2]$ driven by diffusion defined in \eqref{eq:SDEX1} has the property that
  $\tau_\infty<\infty$, a.s. Moreover,
\begin{align}\label{s30.2}
\P[\x](\tau_\infty > t ) \leq c_1\e^{- c_2 t}, \qquad t \geq 0,\ \x \in (0,2]^N,
\end{align}
where $c_1$ and $c_2$ depend only on $N$ and $\beta$,
and satisfy $0< c_1, c_2 < \infty$.
\end{theorem}

\begin{remark}
(i) If we take $\beta = 2$ in \eqref{eq:SDEX1} then the diffusion is a Bessel process (locally near 0). Hence, we may say that Theorem \ref{thm:main2} is concerned with a diffusion with a drift ``slightly stronger'' than the drift of any Bessel process.

(ii) The theorem still holds if the  constant $1/\beta$ in the drift term in \eqref{eq:SDEX1} is replaced by any other positive constant. We chose $1/\beta$ to simplify some formulas in the proof.

(iii) The diffusion \eqref{eq:SDEX1} is reflected at 2 so that we can prove the exponential bound in \eqref{s30.2}. For some Markov processes, the hitting time of a point can be finite almost surely but it may have an infinite expectation; the hitting time of 0 by one-dimensional Brownian motion starting at 1 is a classical example of such situation. The reflection is used in \eqref{eq:SDEX1} to get rid of the effects of excursions of the diffusion far away from the boundary at 0. A different example could be constructed based on a diffusion on $(0,\infty)$ with no reflection but with very strong negative drift far away from 0. 
\end{remark}

We end this section with two open problems.
\begin{problem}
Find necessary and sufficient conditions, in terms of $N$ and $\nu$, for non-extinction in a finite time of an $N$-particle Fleming-Viot process driven by $\nu$-dimensional Bessel process.
\end{problem}

\begin{problem}
Does there exist a Fleming-Viot-type process, not necessarily driven by Brownian motion, such that $\tau_\infty = \infty$, a.s., for the $N$-particle system, but $\tau_\infty < \infty$ with positive probability for the $(N+1)$-particle system, for some $N\geq 2$?
\end{problem}

The rest of the paper contains the proofs of the two main theorems.

\section{Proof of Theorem \ref{thm:main1}}\label{sec:proof1}

\subsection{Bessel processes}\label{sec:bessel}

We start with a review of some facts about Bessel processes and Gamma distributions.
Let $Z_t$, $t\geq 0$, be a square of Bessel process of dimension $\nu\in\R$ starting at
$x\geq 0$, ($Z\sim\besq{x}$, for short), i.e., $Z$ is the unique strong solution to
stochastic differential equation
\begin{equation*}
  dZ_t=\nu \,dt+2\sqrt{|Z_t|}\,dW_t,\quad Z_0=x,
\end{equation*}
where $W$ is a one-dimensional Brownian motion (see \cite[Chapter
11]{RevuzYor99} for the case $\nu\geq 0$ and \cite{GoingYor03} for the general case).

Squares of Bessel processes have the following scaling property: if
$Z_t\sim\besq{x}$ and for some $c>0$ and all $t\geq 0$ we have $Z'_t=cZ_{c^{-1}t}$, then 
$Z'\sim\besq{cx}$.

If $Z\sim\besq{x}$ with $x>0$, and $T_0$ denotes the first hitting time of 0,
then $T_0=\infty$, a.s., if $\nu\geq 2$, and $T_0<\infty$, a.s., if $\nu<2$. Moreover, in the
latter case we have
\begin{equation}\label{eq:invgamma}
  T_0\overset{d}{=}\frac{x}{2G}, 
\end{equation}
where $G$ is $\Gamma\left( 1-\frac{\nu}{2} \right)$-distributed random variable \cite[eqn.~(15)]{GoingYor03}. Here
and in what follows we say that a random variable is $\Gamma(\alpha)$-distributed if it has the density
\begin{equation*}
  f_\alpha(x)=\frac{1}{\Gamma(\alpha)}\,x^{\alpha-1}\e^{-x},\quad x>0,\,\alpha>0,
\end{equation*}
where
\begin{equation*}
  \Gamma(\alpha)=\int_0^\infty x^{\alpha-1}\,\e^{-x}\,dx
\end{equation*}
denotes the standard gamma function. Note that we consider only a one-parameter family of gamma densities, unlike the traditional two-parameter family.

In \cite{RevuzYor99}, Bessel process $X$ of dimension $\nu\geq 0$ starting at $x\geq 0$
($X\sim\bes{x}$), is defined as the square root of $\besq{x^2}$ process $Z$. If $\nu\geq
0$, then by so called comparison theorems, the paths of $Z_t$ are defined for all
$t\geq 0$, so $X_t$ is well defined for all $t\in [0,\infty)$. 
  We define Bessel process $X$ of dimension $\nu<0$ starting at $x\geq 0$ as the square root
  of a $\besq{x^2}$ process $Z$, i.e.,
$X_t=\sqrt{Z_t}$ for $ t\leq T_0$.
For any real $\nu$, these definitions are equivalent to the definition given in \eqref{s30.3} by the It\^o formula.

Processes $\bes{x}$ with $\nu\in \R$ scale as follows. If $X\sim\bes{x}$ is a Bessel process on $[0,T_0)$, then for all
$c>0$, $c X_{c^{-2}t }$ is a $\bes{cx}$ process on $[0,c^2T_0]$. This follows easily from
the scaling property of $\besq{x}$ processes.

\subsection{Proof of Theorem \ref{thm:main1} (i)}

We start with an alternative construction of the Fleming-Viot process $\X$.
Let $X=(X_t,\,t\in[0, T_0))$ be a $\bes{1}$ process. Let $\Y=\br{Y^{1}_t,Y^{2}_t}$, where
$Y^{1}_t$ and $Y^{2}_t$ are independent copies of $X_t$ and let $\Y^i_t =
\br{Y^{i,1}_t,Y^{i,2}_t}$, $i=1,2,\dotsc$, be a sequence of independent copies of $\Y$.
For $i=1,2,\dotsc$ we set
\begin{gather*}
\sigma_i = \inf\setof{t>0}{Y^{i,1}_t \wedge Y^{i,2}_t = 0}, \\
\intertext{and}
\alpha_i = Y^{i,1}_{\sigma_i} \vee Y^{i,2}_{\sigma_i}.
\end{gather*}
It is easily seen that $\sigma_1$ may be represented as $\sigma_1=\min(T_0,T'_0)$, where
$T'_0$ is an independent copy of $T_0$, and that $(\sigma_i,\,i=1,2,\dotsc)$ is a sequence of
independent and identically distributed random variables.

We construct a two-particle Fleming-Viot type process $\X_t = \br{X^1_t, X^2_t}$ as
follows. First let $\tau_1=\sigma_1$ and set $\X_t = \Y^1_t$ for $t\in\halfint{0,\tau_1}$.
At $\tau_1$ one of the particles hits the boundary of $D=(0,\infty)$, and it jumps to
$\xi_1 = \alpha_1$. To continue the process we use the scaling property of $\Y_t$: let
$\tau_2 = \tau_1 + \xi_1^2\sigma_2$ and set $\X_t = \xi_1\Y^2_{\xi_1^{-2}(t-\tau_1)}$
for $t\in\halfint{\tau_1,\tau_2}$. At $\tau_2$, one of the particles hits the boundary and
jumps, this time to $\xi_2 = \alpha_2\xi_1$. We continue the process in the same way by
setting
\begin{gather*}
\xi_j =\prod_{i=1}^j \alpha_i,\\
\tau_n = \sum_{j=1}^{n} \xi_{j-1}^2 \sigma_{j}, \\
\intertext{and}
\X_t = \xi_n \Y^n_{\xi_n^{-2}(t-\tau_n)} ,
\qquad \text{for } t\in\halfint{\tau_n,\tau_{n+1}}.
\end{gather*}
It is easy to see that the construction of $\X$ given above is equivalent to that given in the Introduction, except that the driving process is a $\nu$-dimensional Bessel process. The process $\X_t$ is well defined up until
$\tau_\infty$, and we will show now that $\tau_\infty<\infty$ almost surely  if and only
if $\nu<0$. 

Note that $\X_0 = (1,1)$ for the process constructed above. However,
it is easy to see that for any two starting points $\X_0 = (x_0^1, x_0^2)$
and $\X_0 = (z_0^1, z_0^2)$ with $x_0^1,x_0^2, z_0^1, z_0^2 >0$, 
the distributions of $\X_{\tau_1}$ are mutually absolutely continuous.
This implies that the argument given below proves the theorem
for any initial value of $\X$.

The case $\nu\geq 2$ is very simple: then $\sigma_1=\infty$, a.s.,
so $\tau_\infty=\infty$, a.s. So for the rest of this section we assume that $\nu<2$.

To check whether  $\tau_\infty<\infty$ or $\tau_\infty=\infty$, we will apply the following
theorem. Let
$\log^+ x=\max(\log x, 0)$.

\begin{theorem}\label{thm:diaco_freed}
(\cite{DiaconisFreedman99}; see also \cite{BougerolNico92} or
\cite{GoldieMaller00})
  Let $\left\{ (A_n,B_n),\, n\geq 1 \right\}$ be a sequence of independent and identically
  distributed random vectors such that $A_n,B_n \in \R$ and
  \begin{equation*}
    \E\left(\log^{+}|A_1|\right)<\infty,\quad \E\left(\log^{+}|B_1|\right)<\infty.
  \end{equation*}
  Then the infinite random series
  \begin{equation*}
    \sum_{n=1}^\infty\Bigl( \prod_{j=1}^{n-1} A_j \Bigr) B_n
  \end{equation*}
  converges a.s. to a finite limit if and only if 
  \begin{equation*}
    \E\log|A_1|<0.
  \end{equation*}
\end{theorem}

We will apply Theorem \ref{thm:diaco_freed} with $A_n=\alpha_n^2$ and $B_n=\sigma_n$. Thus, in
order to prove Theorem \ref{thm:main1} (i), it suffices to show that
\begin{enumerate}[(i)]
  \item $\E\log\sigma_1<\infty$ for $\nu<2$;
  \item $\E\log(\alpha_1^2)<0$ for $\nu<0$ and $\E\log(\alpha_1^2)\geq 0$ for $\nu\geq 0$.
\end{enumerate}

\begin{proof}[Proof of (i)]
  Note that, in view of \eqref{eq:invgamma},
  \begin{equation*}
    \E\log\sigma_1  \leq \E\log T_0 =\E\log\frac{1}{2G} =-\log 2- \E\log G,
  \end{equation*}
  where $G\sim\Gamma\left( 1-\frac{\nu}{2} \right)$. But for $G\sim\Gamma(\alpha)$ with $\alpha = 1-\frac{\nu}{2} >0$, we have
  \begin{equation}\label{eq:ElogG}
    \begin{split}
      \E\log G & =\int_0^\infty \log x\, f_\alpha(x)\,dx\\
      & =\frac{1}{\Gamma(\alpha)}\int_0^\infty  x^{\alpha-1}\log x\,\e^{-x}\,dx\\
      & =\frac{1}{\Gamma(\alpha)}\,\frac{d}{d\alpha}\Gamma(\alpha)=\psi(\alpha)<\infty,
    \end{split}
  \end{equation}
  where $\psi$ is well known digamma function
(\cite[Section 6.3]{AS}) defined as
  \begin{equation*}
    \psi(x)=\frac d{dx}\log\Gamma(x)=\frac{\Gamma'(x)}{\Gamma(x)}. \qedhere
  \end{equation*}
\end{proof}

\begin{proof}[Proof of (ii)]
  By Theorem 11 of \cite{Pal10} we get that the density of $\alpha_1^2$ is given by
  \begin{equation*}
    \begin{split}
      h_\nu(y)&=\frac{(y+2)^{\nu-3}}{\Gamma\left( 1-\frac{\nu}{2} \right)}
      \sum_{n=0}^\infty \frac{\Gamma\left( 3-\nu+2n \right)}{n!\Gamma\left(
      2-\frac{\nu}{2}+n \right)}\left( \frac{y}{(y+2)^2} \right)^n\\
      & = \frac{(y+2)^{\nu-3}}{\Gamma\left( 1-\frac{\nu}{2} \right)}
      g\left( \frac{y}{(y+2)^2} \right),
    \end{split}
  \end{equation*}
  where
  \begin{equation*}
    g(z)=\sum_{n=0}^{\infty}c_n z^n
  \end{equation*}
  with
  \begin{equation*}
    c_n=\frac{\Gamma\left( 2n+3-\nu \right)}
    {n!\Gamma\left( n+ 2-\frac{\nu}{2} \right)}, 
    \quad n=0,1,2,\dotsc.
  \end{equation*}

  By the duplication formula for the gamma function (\cite[eqn.~6.1.18]{AS}), i.e.,
  \begin{equation*}
    \Gamma(2z)=\frac{2^{2z-1}}{\sqrt\pi}\Gamma(z)\Gamma\left( z+\frac{1}{2} \right),
  \end{equation*}
  we have 
  \begin{equation*}
    \begin{split}
      c_n & = \frac{2^{2n+2-\nu}}{\sqrt\pi}\cdot\frac{\Gamma\left( n+\frac{3-\nu}{2}
      \right)}{n!}\\ 
      & =\frac{2^{2n+2-\nu}}{\sqrt\pi}\binom{n+\frac{1-\nu}{2}}n \Gamma\left(
      \frac{3-\nu}{2} \right),
     \end{split}
  \end{equation*}
  where for $x>-1$ and $k\in\N$,
  \begin{equation*}
    \binom x k=\frac{\Gamma(x+1)}{k!\Gamma(x-k+1)}
  \end{equation*}
  is a generalized binomial coefficient. Therefore, 
  \begin{equation*}
    \begin{split}
      g(z) & = \frac{2^{2-\nu}}{\sqrt\pi}\Gamma\left( \frac{3-\nu}{2} \right)
      \sum_{n=0}^\infty \binom{n+\frac{1-\nu}{2}}n (4z)^n\\ 
       & = \frac{2^{2-\nu}}{\sqrt\pi}\Gamma\left( \frac{3-\nu}{2} \right)
       \left( \frac{1}{1-4z} \right)^{\frac{3-\nu}{2}},
     \end{split}
  \end{equation*}
  as for $a\in\R$
  \begin{equation*}
    \sum_{n=0}^\infty\binom{n+a}n z^n= (1-z)^{-a-1}.
  \end{equation*}
  Now
  \begin{equation*}
    g\left( \frac{y}{(y+2)^2} \right)=\frac{2^{2-\nu}}{\sqrt\pi}\Gamma\left( \frac{3-\nu}{2} \right)
    \frac{(y+2)^{3-\nu}}{(y^2+4)^{\frac{3-\nu}{2}}}
  \end{equation*}
  and therefore for $y\geq 0$
  \begin{equation*}
    h_\nu(y)=\frac{2^{2-\nu}}{\sqrt\pi}\frac{\Gamma\left( \frac{3-\nu}{2}
    \right)}{\Gamma\left( 1-\frac{\nu}{2} \right)}\,\frac{1}{\left( y^2+4
    \right)^{\frac{3-\nu}{2}}}.
  \end{equation*}
  So, to prove (ii) we need to study the sign of the integral
  \begin{equation*}
    I(\nu)=\int_0^{\infty}h_\nu(y)\log y\,dy.
  \end{equation*}
Recall the Student's $t$-distribution with $a>0$ degrees of
  freedom (\cite[section 26.7]{AS}).
  The density for this distribution is given by
  \begin{equation*}
    f(x;a)=
    \frac{\Gamma\left( \frac{a+1}{2} \right)}{\sqrt{\pi a}\,\Gamma\left(
    \frac{a}{2} \right)}\left( 1+\frac{x^2}{a} \right)^{-\frac{a+1}{2}},\quad -\infty < x < \infty.
  \end{equation*}
  Changing the variable $y=\frac{2x}{\sqrt{2-\nu}}$ in $I(\nu)$ we get
  \begin{equation*}
    \begin{split}
       I(\nu)&=\int_0^{\infty}f(x;2-\nu)\log\frac{2x}{\sqrt{2-\nu}}\,dx\\
      &=\frac{1}{2}\E\log\frac{2\abs{X}}{\sqrt{2-\nu}},
    \end{split}
  \end{equation*}
  where $X$ is a random variable with $t$-distribution with $(2-\nu)$-degrees of freedom.
  
  It is well known (\cite[section 26.7]{AS}) that
  \begin{equation*}
    X\overset{d}{=}\frac{Z\sqrt{2-\nu}}{\sqrt{V}},
  \end{equation*}
  where $Z$ has standard normal distribution and $V$ has chi-squared distribution with
  $(2-\nu)$ degrees of freedom, and $Z$ and $V$ are independent. Therefore
  \begin{equation*}
    \begin{split}
      I(\nu)&=\frac12\E\log\frac{2\abs{Z}}{\sqrt{V}}\\
      &=\frac12\log 2+\frac{1}{4}\left( \E\log\frac{Z^2}{2}-\E\log\frac{V}{2} \right).
    \end{split}
  \end{equation*}
  Note that $\frac{Z^2}{2}$ has $\Gamma\left(\frac{1}{2}\right)$ distribution and
  $\frac{V}{2}$ has $\Gamma\left(\frac{2-\nu}{2} \right)$ distribution. Therefore, by
  \eqref{eq:ElogG},
  \begin{equation*}
      I(\nu)=\frac12\log 2+\frac{1}{4}\left( \psi\left( \frac{1}{2} \right)
    -\psi\left( \frac{2-\nu}{2} \right) \right).
  \end{equation*}
The function $\psi$ is strictly increasing with $\psi\left( \frac{1}{2} \right)=-2\log 2-\gamma$
  and $\psi(1)=-\gamma$ where $\gamma$ is the Euler constant (\cite[eqns.~6.3.2, 6.3.3]{AS}). Using these facts we see that
  \begin{equation*}
       I(\nu)=\frac{1}{4}\left( \psi(1)-\psi\left( \frac{2-\nu}{2} \right) \right),
  \end{equation*}
  and therefore $I(\nu)<0$ iff $\frac{2-\nu}{2}>1$ iff $\nu<0$. This completes the proof
  of (ii) and of Theorem \ref{thm:main1} (i).
\end{proof}

\subsection{Proof of Theorem \ref{thm:main1} (ii)}

Suppose that $\X= (X^1, \dots , X^N)$ is a Fleming-Viot process driven by $\nu$-dimensional Bessel process, $\X_0 = (x^1, \dots, x^N)$, $x^j >0$ for all $1\leq j \leq N$, and $N\nu \geq 2$. Let $Z_t = (X^1_t)^2 + \cdots + (X^N_t)^2$ and $z_0 = (x^1)^2 + \cdots + (x^N)^2>0$.
According to \cite[Thm. 2.1]{Shiga_Wat}, the process $\{Z_t, t\in[0, \tau_1)\}$
is an $(N\nu)$-dimensional square of Bessel process, i.e., it has distribution $\besselsq^{N\nu}(z_0)$. More generally, 
$\{Z_t, t\in[\tau_k, \tau_{k+1})\}$ has distribution $\int\besselsq^{N\nu}(z)\P(Z_{\tau_k} \in dz)$ for $k\geq 0$, where,
by convention, $\tau_0=0$.
Let $Y_t = Z_t^{1/2}$, $B_0=0$ and define $B$ inductively on intervals 
$(\tau_k, \tau_{k+1}]$ 
by  
\begin{align*}
B_t = B_{\tau_k} + \int_{\tau_k}^t dY_s - \int_{\tau_k}^t\frac{N\nu-1}{2 Y_s} ds.
\end{align*}
Then, by the It\^o formula,
$B$ is a Brownian motion and 
\begin{equation*}
  dZ_t=2\sqrt{Z_t}\,dB_t + N\nu \,dt,
\end{equation*}
for $t\in (\tau_k, \tau_{k+1})$, $k\geq 0$.
Let $\wh Z^t$ be defined by $\wh Z_0 = z_0$ and
\begin{equation*}
  \wh Z_t= \int_0^t 2\sqrt{\wh Z_s}\,dB_s + N\nu t, \qquad t\geq 0.
\end{equation*}
By definition, $\wh Z$ is an $(N\nu)$-dimensional squared Bessel process on $[0, \infty)$. We assumed that $N\nu \geq 2$ and $z_0>0$ so we have $\wh Z_t >0$ for all $t\geq 0$, a.s.
Since $\wh Z$ is continuous, for every integer $j>0$ there exists a random variable $a_j$ such that $\wh Z_t > a_j >0$ for all $t\in[0,j]$, a.s.
Note that $\wh Z_t = Z_t$ for $t\in [0, \tau_1)$ and $\wh Z_{\tau_1} < Z_{\tau_1}$ because $Z$ has a positive jump at time $\tau_1$. Strong existence and uniqueness for SDE's with smooth coefficients implies that $\wh Z_t \leq Z_t$ for all $t\in [\tau_1, \tau_2)$, because if the trajectories of $\wh Z$ and $Z$ ever meet then they have to be identical after that time up to $\tau_2$. Once again, 
$\wh Z_{\tau_2} < Z_{\tau_2}$ because $Z$ has a positive jump at time $\tau_2$. By induction, $\wh Z_t \leq Z_t$ for all $t\in [\tau_k, \tau_{k+1})$, $k\geq 0$, a.s. 
Hence, $ Z_t > a_j >0$ for all $t\in[0,j]$ and $j>0$, a.s. This implies that $\tau _\infty = \infty$, a.s., by an argument similar to that in Lemma 5.2 of 
\cite{BieniekBurdzyFinch10}.

\section{Proof of Theorem \ref{thm:main2}}

\subsection{Preliminaries}

We will give new meanings to some symbols used in the previous section. Constants denoted by $c$ with subscripts will be tacitly assumed to be strictly positive and finite; in addition, they may be assumed to satisfy some other conditions.

(i) Let $W$ be one-dimensional Brownian motion  and let $b$ be
    a Lipschitz function defined on an interval in $\R$, i.e., $|b(x_1) - b(x_2)| \leq L |x_1 - x_2|$ for some $L< \infty$ and all $x_1$ and $x_2$ in the domain of $b$.  Consider a diffusion $X_t, t\in [s,u]$, satisfying the following stochastic differential equation,
    \begin{equation}\label{o6.1}
      dX_t= dW_t+b\left( X_t \right)\,dt,\quad
      X_s=a.
    \end{equation}
Let $y_t$ be the solution to the ordinary differential equation
    \begin{equation*}
      \frac d {dt} y_t=b(y_t),\quad y_s=a.
    \end{equation*}
We will later write
$y'=b(y)$ instead of $\frac d {dt} y_t=b(y_t)$.

The following inequality appears in Ch.~3, Sect.~1 of the book by Freidlin and Wentzell \cite{FredlinWentzell83}. For every $\delta>0$,
    \begin{align*}
      \P\left( \sup_{s\leq t\leq u}\left|X_t-y_t\right|>\delta \right)\leq 
\P \left( \sup_{s\leq t\leq u} |W_t| > \delta \e^{-L(u-s)}
\right),
    \end{align*}
where $L$ is a Lipschitz constant of $b$. It follows that
    \begin{align}\nonumber
      \P\left( \sup_{s\leq t\leq u}\left|X_t-y_t\right|>\delta \right)
&\leq 
\P \left( \sup_{s\leq t\leq u} W_t > \delta \e^{-L(u-s)}
\right)
+\P \left( \inf_{s\leq t\leq u} W_t <- \delta \e^{-L(u-s)}
\right)\\
& = 2\P \left( \sup_{s\leq t\leq u} W_t > \delta \e^{-L(u-s)}
\right)\nonumber\\
& = 4\P \left(  W_u-W_s > \delta \e^{-L(u-s)}
\right)\nonumber\\
&\leq 
      c_0 \exp\left( -\frac{\delta^2}{2(u-s)}\,\e^{-2L(u-s)} \right), \label{eq:FWineq}
    \end{align}
where $c_0$ is an absolute constant.

(ii) Recall that $\beta > 2$ and consider the function
    \begin{equation}\label{eq:defb}
      b(x)=-\frac{1}{\beta x^{\beta-1}},\quad x>0.
    \end{equation}
We need the assumption that $\beta >2$ for the main part of the argument but
many calculations given below hold for a larger family of $\beta$'s. It is easy to check that 
    \begin{equation}\label{eq:yt}
      y_{s,a}(t):=\left( a^\beta+s -t \right)^{1/\beta},\quad s\leq t\leq
      s+a^\beta,
    \end{equation}
is the solution to the ordinary differential equation
    \begin{equation}\label{eq:ODEb}
      y'=b(y)
    \end{equation}
    with the initial condition 
    $y_{s,a}(s)=a$, 
    where $s\in \R$, $a>0$.
Note that the function $y_{s,a}(t)$ approaches 0 vertically at $ t = s+a^\beta$. 

(iii) 
Fix any
    $\gamma\in(0,1)$ and let $L$ be the Lipschitz constant of $b$ on the interval
    $\left[ a(\gamma/2)^{1/\beta}/2, 2a\right]$. Then
    \begin{equation*}
L=b'\bigl(  a(\gamma/2)^{1/\beta}/2 \bigr), \qquad
      b'(x)=\frac{\beta-1}{\beta x^\beta},
\end{equation*}
and, therefore,
\begin{equation}\label{o7.1}
   L=\frac{\beta-1}{\beta(a(\gamma/2)^{1/\beta}/2)^\beta}
=
\frac{\beta-1}{\beta\gamma 2^{1-\beta} a ^\beta}.
    \end{equation}
Let $X$ be the solution to \eqref{o6.1} with $b$ defined in \eqref{eq:defb}.
Assume that $\delta>0$ is so small that 
\begin{align}\label{o6.5}
a(\gamma/2)^{1/\beta}/2 \leq
y_{0,a}((1-\gamma/2)a^\beta)-\delta < y_{0,a}(0) + \delta \leq 2 a.
\end{align}
It follows 
that if $ \sup_{0\leq t\leq (1-\gamma/2)a^\beta} \left|X_t-y_{0,a}(t)\right|\leq\delta$
then $X_t\in  \left[ a(\gamma/2)^{1/\beta}/2, 2a\right]$ 
for $0\leq t\leq (1-\gamma/2)a^\beta$.
Hence, we can apply \eqref{eq:FWineq} with $L$ given by \eqref{o7.1} to obtain the following estimate
    \begin{equation}\label{eq:y0A}
     \P\left( \sup_{0\leq t\leq (1-\gamma/2)a^\beta} \left|X_t-y_{0,a}(t)\right|>\delta 
           \mid X_0 =a \right)
       \leq c_0 \exp\left( -c_1\,\frac{\delta^2}{a^\beta} \right),
    \end{equation}
    where 
    $X_t$ satisfies \eqref{o6.1} and   
    $c_1$ depends on $\beta$ and $\gamma$, but it does not depend on $\delta$ and
    $a$.

(iv) 
Suppose that $a>0$ and $u=(1-\gamma)a^\beta$. Then
    $y_{0,a}(u)=a\gamma^{1/\beta}$. 
    For $\varepsilon\in(0,1)$, let 
    \begin{equation}\label{o1.1}
\bar\delta=      \bar\delta(\varepsilon)=a\left[ 1-(1-\varepsilon\gamma)^{1/\beta} \right].
    \end{equation}
We fix $\eps\in(0,1)$ so small that for all $a>0$ the inequality \eqref{o6.5} is satisfied with $\bar\delta$ in place of $\delta$ and, moreover,
\begin{align}\label{o4.1}
(1-\gamma/2)(a-\bar\delta)^\beta >(1-\gamma)a^\beta
\end{align}
and
\begin{align}\label{o4.7}
\gamma^{1/\beta}  +  1-(1-\varepsilon\gamma)^{1/\beta} < 1.
\end{align}
If $\delta\in[-\bar\delta,\bar\delta]$ then
    \begin{equation*}
      y_{0, a+\delta}(t)=\left( \left( a+\delta \right)^{\beta}-t
      \right)^{1/\beta},\quad 0\leq t\leq (a+\delta)^{\beta}.
    \end{equation*}
It is straightforward to check that 
    \begin{equation*}
      y_{0, a+\delta}(u)\geq a\gamma^{1/\beta}(1-\varepsilon)^{1/\beta}>0.
    \end{equation*}
We will estimate the difference $ y_{0, a+\delta}(u) - y_{0, a}(u)$ as a function of $\delta$.
    Let 
    \begin{equation*}
      f(\delta)=y_{0, a+\delta}(u)
      =\left( \left( a+\delta \right)^{\beta}-(1-\gamma)a^\beta
      \right)^{1/\beta}.
    \end{equation*}
    Then
    \begin{equation*}
y_{0, a+\delta}(u) - y_{0, a}(u) =f(\delta)-f(0)=\delta f'(\hat\delta),
    \end{equation*}
    where $\hat\delta$ is between 0 and $\delta$. But
    \begin{equation*}
      f'(\delta)=\left( a+\delta \right)^{\beta-1}
        \left( \left( a+\delta \right)^{\beta}-(1-\gamma)a^\beta
        \right)^{1/\beta-1}    
    \end{equation*}
    and
    \begin{equation*}
      f''(\delta)=-a^\beta(\beta-1)(1-\gamma)\left(
      a+\delta \right)^{\beta-2}\left( \left( a+\delta
      \right)^{\beta}-(1-\gamma)a^\beta \right)^{1/\beta-2}.
    \end{equation*}
It follows from \eqref{o1.1} and the fact that $\eps\in(0,1)$ that $\left( a+\delta \right)^{\beta}-(1-\gamma)a^\beta>0$.
Thus $f''<0$ and $f'$ is strictly decreasing. Therefore, for
    $\delta\in[-\bar\delta,\bar\delta]$ we have
    \begin{equation}\label{eq:Mdef}
      f'(\delta)\leq f'(-\bar\delta)=\left(
      \frac{1-\varepsilon\gamma}{\gamma(1-\varepsilon)} \right)^{1-1/\beta}=:M >1.
    \end{equation}
Hence $      \left| y_{0, a+\delta}(u) - y_{0, a}(u) \right|\leq M|\delta|$ for $\delta\in[-\bar\delta,\bar\delta]$.
Suppose that $\delta'\in(0,\bar\delta]$ and
    $\delta=\frac{\delta'}{M+1}$. If $a_0\in[a-\delta,a+\delta]$, then 
    \begin{equation}\label{eq:Mdelta}
      y_{0,a}(u)-\delta' \leq y_{0,a_0}(u) \leq y_{0,a}(u)+\delta'.
    \end{equation}
The function $b(x)$ is strictly increasing on $(0,\infty)$. This easily implies that if $0 < a_1 < a_2$ and $0 \leq s < t \leq a_1^\beta$ then 
\begin{align}\label{o4.3}
y_{0,a_2}(s) - y_{0,a_1}(s) < y_{0,a_2}(t) - y_{0,a_1}(t).
\end{align}
Hence, inequality \eqref{eq:Mdelta} holds in fact for all $t\in[0,u]$ in place of $u$ and, moreover, $y_{0,a-\bar\delta}(t)\leq y_{0,a}(t)-\delta'$. 
Another consequence of \eqref{eq:Mdelta} and \eqref{o4.3} is that
if $s\in[0,u]$ and $a_0\in[y_{0,a}(s)-\delta,y_{0,a}(s)+\delta]$ then for $t \in[s,u]$,
    \begin{equation}\label{o3.5}
      y_{0,a-\bar\delta}(t)\leq
      y_{0,a}(t)-\delta' \leq y_{s,a_0}(t) \leq y_{0,a}(t)+\delta'.
    \end{equation}
The following generalization of \eqref{o4.1} also follows from \eqref{o4.3}. If $s\in[0,u]$ and $a_0 \geq y_{0,a}(s)-\delta$ then
\begin{align}\label{o4.2}
s+(1-\gamma/2)a_0^\beta >(1-\gamma)a^\beta.
\end{align}

\subsection{Proof of Theorem \ref{thm:main2}}

Suppose that $\X_t=(X_t^1,\dotsc,X_t^N)$
is a
Fleming-Viot process on $(0,2]$ driven by the diffusion defined in \eqref{eq:SDEX1}, with an arbitrary $2 \leq N < \infty$.
Recall that the role of the boundary is played by the point 0, and only this point.
In other words, the particles jump only when they approach 0.
    
{\it Step 1}.
 Let $\x=\left( x^1,\dotsc,x^N \right)$, $[N]=\set{1,\dotsc,N}$, and let $j_1$ be the smallest integer in $[N]$ with
\begin{equation*}
      x^{j_1}=\max_{1\leq j\leq N} x^j.
    \end{equation*}
Consider process $\X$ starting from $\X_0=\x$ and let $I_1=\set{j_1}$, $J_1=[N]\setminus I_1$ and $S_0=0$. Let $u=(1-\gamma)(x^{j_1})^\beta$ and 
    \begin{equation*}
      S_1= u \wedge\inf\set{t\geq 0:\exists_{j\in J_1}\,X_t^j=X_t^{j_1}}.
    \end{equation*}
Note that two processes $X^i$ and $X^j$ can meet either when their paths intersect at a time when both processes are continuous or when one of the processes jumps onto the other. Let $j_2$ be the smallest index in $J_1$ such that the equality in the definition of $S_1$ holds with $j = j_2$.
Let $I_2=\set{j_1,j_2}$ and $J_2=[N]\setminus I_2$. 

    Next we proceed by induction. Assume that, for some $n<N$, the sets $I_1,\dotsc,I_n$, $J_1,\dotsc,J_n$,
    and stopping times $S_1<S_2<\dotsc<S_{n-1}$ are defined. Then we let
    \begin{equation*}
      S_{n}= u \wedge\inf\set{t\geq S_{n-1}:\exists_{i\in I_n}\,\exists_{j\in J_n}\,
      X_t^i=X_t^j},
    \end{equation*}
$I_{n+1}=I_n\cup\set{j_{n+1}}$ and $J_{n+1}=[N]\setminus I_{n+1}$, where $j_{n+1}$ is the smallest index in $J_n$ such that the equality in the definition of $S_n$ holds with $j = j_{n+1}$.

The set $I_n$ has $n$ elements which are indices of particles which are ``descendants'' of the particle $X^{j_1}$ that was the highest at time 0. By convention, we let $I_n=I_N$ and
    $S_n=u$ for $n\geq N$.

{\it Step 2}. Write $a=x^{j_1}$ and  $u=(1-\gamma)a^\beta$. Then $\x\in(0,a]^N$. 
Recall $\bar\delta$ and $M$ defined in \eqref{o1.1} and \eqref{eq:Mdef}, and for $1\leq n\leq N$ define 
    \begin{equation}\label{eq:hatdelta}
      \hat\delta_{n}=\frac{\bar\delta}{(M+1)^{N-n}}.
    \end{equation}
Note that $\hat\delta_{n}=(M+1)\hat\delta_{n-1}$.
Consider events
    \begin{equation*}
      \eF_n=\bigcup_{ j\in I_n}
\left\{
\sup_{S_{n-1}\leq
      t<S_n}\left|X_t^j-y_{0,a}(t)\right|>\hat\delta_n 
\right\} .
    \end{equation*}
Note that for every $t$, $\max_{j\in I_n} X^j_t \geq 
\max_{j\in J_n} X^j_t$. Hence,  
\begin{align*}
\bigcup_{1\leq j\leq N} 
\left\{\sup_{0\leq t<u}X_t^j- y_{0,a}(t) >\bar\delta\right\}
\subset  \bigcup_{1\leq n\leq N} \eF _n,
\end{align*}
and, therefore,
    \begin{equation}\label{eq:D}
      \begin{split}
	 \P[\x]&\left( \bigcup_{1\leq j\leq N} 
	 \left\{\sup_{0\leq t<u}X_t^j- y_{0,a}(t) >\bar\delta\right\}
	 \right)\\
&\leq \P[\x] \left( \bigcup_{1\leq n\leq N} \eF _n \right)
\\
&= \P[\x] \left( \bigcup_{1\leq n\leq N} \eF _n \cap \eF_1^c \cap \dots
\cap \eF_{n-1}^c \right)\\
&\leq
\sum_{1\leq n\leq N}
\P[\x] \left(  \eF _n \cap \eF_1^c \cap \dots
\cap \eF_{n-1}^c \right)
\\
	&\leq \sum_{1\leq n\leq N} \P[\x]\left( \eF_n \mid
	\eF_{1}^c,\dotsc,\eF_{n-1}^c \right)\\
	&\leq\sum_{1\leq n\leq N}\sum_{j\in I_n}\P[\x]\left( 
	 \sup_{S_{n-1}\leq t<S_n}
	 \left|X_t^j-y_{0,a}(t)\right|>\hat\delta_n
	 \Bigm|	\eF_{1}^c,\dotsc,\eF_{n-1}^c \right),
      \end{split}
    \end{equation}
    where we adopted the convention
    $\P[\x](\eF_1\mid\eF_0^c)=\P[\x](\eF_1)$.
    
Suppose that $\eF_{n-1}^c$ holds and $j \in I_{n}$. Then $|X^j_{S_{n-1}} - y_{0,a}(S_{n-1})| \leq \hat\delta_{n-1}$.
Let $y^j_t$, $t\geq S_{n-1}$, be a
    solution to $y'=b(y)$ with $y^j_{S_{n-1}}=X^j_{S_{n-1}}$. By \eqref{o3.5}, 
    \begin{equation*}
      \begin{split}
	\P[\x]\left( \sup_{S_{n-1}\leq t<S_{n}} \left| X_t^j- y_{0,a}(t)\right|>\hat\delta_{n} \mid \eF_{n-1}^c
	\right)
	  &\leq \P[\x]\left( \sup_{S_{n-1}\leq t<S_{n}} \left| X_t^j-y_t^j\right|>\hat\delta_{n}
 \mid \eF_{n-1}^c	 \right).
      \end{split}
    \end{equation*}
It follows from \eqref{o4.2} that we can apply \eqref{eq:y0A} (with an appropriate shift of the time scale) to $X^j$, assuming that $\eF_{n-1}^c$ holds, on the interval $[S_{n-1},S_{n}]\subset [0, u] = [0,(1-\gamma)a^\beta]$. 
We obtain
\begin{align*}
&\P[\x]\left( \sup_{S_{n-1}\leq t<S_{n}} \left| X_t^j-y_t^j\right|>\hat\delta_{n}
 \mid \eF_{n-1}^c	 \right)
\leq c_0 \exp\left( -c_2\,\frac{\hat\delta_{n}^2}
{(y_{0,a}(S_{n-1}) + \hat\delta_{n-1})^\beta} \right)\\
&\quad \leq c_0\exp\left( -c_2\,\frac{\bar\delta^2 (M+1)^{-2N}}
{(a + \bar\delta)^\beta} \right)
= c_0 \exp\left( -c_3\,\frac{\bar\delta^2 }
{(a + \bar\delta)^\beta} \right).
\end{align*}
We combine this estimate with \eqref{eq:D} to see that
    \begin{align}\label{o6.2}
	 \P[\x]&\left( \bigcup_{1\leq j\leq N} 
	 \left\{\sup_{0\leq t<u}X_t^j- y_{0,a}(t) >\bar\delta\right\}
	 \right)\\
	&\leq\sum_{1\leq n\leq N}\sum_{j\in I_n}\P[\x]\left( 
	 \sup_{S_{n-1}\leq t<S_n}
	 \left|X_t^j- y_{0,a}(t) \right|>\hat\delta_n
	 \Bigm|	\eF_{1}^c,\dotsc,\eF_{n-1}^c \right)
\nonumber \\
&\leq c_0 N^2 \exp\left( -c_3\,\frac{\bar\delta^2 }
{(a + \bar\delta)^\beta} \right). \nonumber
    \end{align}

{\it Step 3}. We will prove that there exist $v<\infty$ and $r\in(0,2)$ such that if
    $\x\in(0,r]^N$ then
    \begin{equation}\label{o2.1}
      \P[\x](\tau_\infty> v)\leq \frac{1}{2}.
    \end{equation}

Consider an $r\in(0,2)$ and for $\x = (x^1, \dots, x^N)$, 
let $A_0 = \max_j x^j$, $U_0=0$,  and for $k=0,1,2,\dotsc$,
    let
    \begin{align*} 
       U_{k+1}&= U_k+(1-\gamma)  A_k^\beta,\\
  A_{k+1}&= \max _{1\leq j \leq N} X^j_{ U_{k+1}}.
    \end{align*} 
Let $\yk_t$ denote the solution to ODE \eqref{eq:ODEb} with the initial condition
	$\yk_{ U_k}=  A_k$. 
Recall $\eps$ from \eqref{o4.1} and
let $ \Delta_k =  A_k\left[ 1-(1-\varepsilon\gamma)^{1/\beta} \right]$. For $k=0,1,2,\dotsc$ define events
    \begin{equation*}
      \Gamma_k=\left[\max_{1\leq j\leq N} 
      \sup_{ U_{k}\leq t< U_{k+1}} X_t^j>\yk_t+ \Delta_k\right].
    \end{equation*}

Note that $\yk_{U_{k+1}} = \gamma^{1/\beta} A_k$. 
Suppose that 
$\bigcap_{k=0}^\infty\Gamma_k^c$ holds. Then 
\begin{align*}
  A_{k+1} \leq \gamma^{1/\beta}   A_k +   A_k\left[ 1-(1-\varepsilon\gamma)^{1/\beta} \right]
=  c_4  A_k, 
\end{align*}
for all $k$, where $c_4 = \gamma^{1/\beta}  +  1-(1-\varepsilon\gamma)^{1/\beta} < 1$, by \eqref{o4.7}. Hence, $  A_k \leq c_4^k   A_0$ and, therefore, $\sum_k   A_k^\beta < \infty$. 
If we let $v= (1-\gamma) \sum_{k=1}^\infty r ^\beta c_4^{k\beta} < \infty$
then $ \lim _{k\to \infty}  U_k
\leq (1-\gamma) \sum_{k=1}^\infty A_0 ^\beta c_4^{k\beta} \leq v $ and 
$\limsup_{t \uparrow v}
\max_{1\leq j\leq N}  X_t^j = 0$. This implies easily that
$\tau_\infty\leq v$.
Thus, to prove \eqref{o2.1}, it will suffice to show that
there exists $r\in(0,2)$ such that if
    $\x\in(0,r]^N$, then
    \begin{equation}\label{eq:mainseries}
      \P[\x]\left( \bigcup_{k=0}^\infty\Gamma_k \right)<\frac{1}{2}.
    \end{equation}
    But
    \begin{equation}\label{eq:A}
      \P[\x]\left( \bigcup_{k=0}^\infty\Gamma_k \right)\leq 
      \P[\x](\Gamma_0)+\sum_{k=1}^\infty \P[\x]\left( \Gamma_k \mid
      \Gamma_0^c,\dotsc,\Gamma_{k-1}^c \right).
    \end{equation}
By \eqref{o6.2} and the strong Markov property applied at $U_k$,
    \begin{equation*}
      \begin{split}
	\P[\x]&\left( \Gamma_k \mid \Gamma_0^c,\dotsc,\Gamma_{k-1}^c \right)
	\leq c_0 N^2\exp\left( -c_3 \frac{ \Delta_k^2}{( A_k +   \Delta_k)^\beta} \right) \\ 
	& = c_0 N^2\exp\left( - c_3 A_k^{2-\beta}
	\frac{(1-(1-\varepsilon\gamma)^{1/\beta})^2}{(2-(1-\varepsilon\gamma)^{1/\beta})^\beta} \right)\\
	& \leq c_0 N^2\exp\left( - c_5  A_0^{2-\beta} c_4^{k(2-\beta)} \right),
      \end{split}
    \end{equation*}
    where $c_4 < 1$. So by \eqref{eq:A}, if $\max_j x^j\leq r$, then 
    \begin{equation*}
	\P[\x]\left( \bigcup_{k=0}^\infty \Gamma_k \right)
	\leq c_0 N^2\sum_{k=0}^\infty\exp\left( - c_5   r^{2-\beta} c_4^{k(2-\beta)} \right),
    \end{equation*}
    which is convergent.  Since $2-\beta<0$, we can
    choose $r>0$ so small that the above sum is less than $1/2$, proving
    \eqref{eq:mainseries}. 

{\it Step 4}. Let $r\in(0,2)$ and $v$ be as in Step 3. Partition the set $(0,2]^N$ into two sets $A=(0,r]^N$ and $A^c$. First we will show that the time when process $\X$ enters the set $A$ has a distribution with an exponentially decreasing tail.

    So assume that $\X_0\in A^c$ and let 
    \begin{equation*}
      I_1=\set{j\in[N]:X_0^j\in(r,2]},\quad I_2=[N]\setminus I_1.
    \end{equation*}
    Let $\tau_1^{j}$ be the the first hitting time of $0$ by the process $X^j$ and let
    \begin{equation*}
      \eta=
      \begin{cases}
	0 & \text{if $I_2=\emptyset$,}\\
	\max_{j\in I_2}\left\{ \tau_1^j \right\},&\text{otherwise}.
      \end{cases}
    \end{equation*}
    Consider
    \begin{equation*}
	p_1(\x)=\P[\x]\left\{ \forall_{j\in I_1}\,\forall_{0\leq t\leq 1/2}\,
	X_t^j\in\left[ \frac{r}{2},2 \right];\eta<1/2;
	\forall_{i\in I_2}\forall_{\tau_1^j<t\leq 1/2}\, 
	X_t^j\in\left[ \frac{r}{4},2 \right]\right\}.
    \end{equation*}
We will argue that for $\x\in A^c$ we have 
    \begin{equation}\label{eq:p1}
      p_1(\x)\geq p_1>0. 
    \end{equation}
    Indeed, with probability at least $q_1>0$ any  particle from $I_1$ stays in the
    interval $[r/2,2]$ up to time $t=1/2$. With probability at least $q_2>0$ any particle
    from $I_2$ hits $0$ before time $t=1/2$; with probability at least $1/N$ it
    jumps onto a particle in $I_1$; and then with probability at least
    $q_3>0$ it stays in the interval $[r/4,2]$ up to time $t=1/2$. Therefore \eqref{eq:p1}
    holds with $p_1=(q_1q_2q_3/N)^N$. 
Obviously $q_1, q_2, q_3$ and
    $p_1$ depend on $r$.

    Next, if we define
    \begin{equation*}
      p_2(\x)=\P[\x]\left\{ \forall_{j\in[N]}\,\forall_{0<t<1/2} \, 
      X_t^j\in\left[ \frac{r}{8},2 \right];
      X_{1/2}^j\in\left[ \frac{r}{8},r \right]\right\},
    \end{equation*}
    and $B=\left[ \frac{r}{4},2 \right]^N$, then an argument similar to that proving \eqref{eq:p1} shows
    that for $\x\in B$ we have $p_2(\x)\geq p_2>0$, where  $p_2$ depends on $r$.
    Therefore, by the Markov property at time $t=1/2$, for $\x\in A^c$ we have 
    \begin{equation}\label{o6.4}
      \P[\x](\X_1\in A)\geq p:=p_1p_2>0.
    \end{equation}
    Now let
    \begin{equation*}
      T=\inf\left\{ t\geq 0:\X_{ t}\in A \right\}.
    \end{equation*}
By \eqref{o6.4}, for all $\x\in(0,2]^N$, 
    \begin{equation*}
      \P[\x](T\leq 1)\geq p>0.
    \end{equation*}
    Applying the Markov property at $t=1,2, \dots$ we obtain
    \begin{equation*}
	\P[\x](T\geq k) \leq (1-p)^{k}.
    \end{equation*}

 Choose $k$ so large that $(1-p)^{k}<\frac{1}{2}$. Recall that $r$ and $v$ are as in Step 3. Let $\theta$ denote the usual Markovian shift operator. Then for any $\x\in (0,2]^N$,
    \begin{equation*}
      \P[\x]\left( \tau_\infty\geq k+v \right)\leq\P[\x](T\geq
      k)+\P[\x](\tau_\infty \circ \theta_T \geq v)\leq (1-p)^{k}+\frac{1}{2}:=q<1.
    \end{equation*}
    Therefore, applying the Markov property at times $k+v, 2(k+v), 3(k+v), \dots$, we obtain,
    \begin{equation*}
      \P[\x](\tau_\infty \geq n(k+v))\leq q^n,\quad n=1,2,\dotsc,
    \end{equation*}
which proves \eqref{s30.2}. This implies that $\tau_\infty < \infty$, a.s.

\bibliographystyle{abbrv}
\bibliography{bessel}

\end{document}